\documentclass[reprint,aps,prl,superscriptaddress,amsmath,amssymb]{revtex4-2}

\usepackage{graphicx}
\usepackage{bm}
\usepackage{hyperref}

\begin{document}

\title{Hidden Harmonic Structure, Universal Damping, and Stability Bounds in Nonlinear Contact Dynamics}

\author{Y. T. Feng}
\email{y.feng@swansea.ac.uk}
\affiliation{Zienkiewicz Institute for Modelling, Data and AI\\
	Faculty of Science and Engineering, Swansea University, Swansea SA1 8EN, UK}

\date{\today}

\begin{abstract}
	Nonlinear contact dynamics are widely regarded as intrinsically nonlinear systems whose behaviour depends strongly on geometry and impact conditions. Here we show that any one-dimensional conservative contact system satisfying monotone energy-consistent conditions admits two complementary structures: (i) a canonical action–angle representation in physical time, and (ii) an exact harmonic oscillator representation under an energy-based coordinate transformation combined with time reparametrisation. This reveals a hidden linear structure underlying nonlinear contact interactions. Building on this result, we derive a unique universal damping law that preserves linear dissipative dynamics in the transformed harmonic space, and establish a rigorous, closed-form lower bound for the critical timestep in numerical simulations. The framework generalises classical power-law contact models and provides a unified basis for restitution control across arbitrary geometries, recovering known exact solutions as explicit monomial special cases.
\end{abstract}

\maketitle

\section{1. Introduction}

Contact interactions are the fundamental mechanism of momentum and energy transfer across a vast range of physical systems. While central to macroscopic granular mechanics and impact dynamics \cite{hertz1882, tsuji1992}, the same nonlinear contact laws govern solitary wave propagation in acoustic metamaterials \cite{daraio2006}, the biomechanical probing of soft matter via atomic force microscopy (AFM) \cite{giessibl2003}, and nanoscale tribology \cite{garcia2002}. Despite their ubiquity, these nonlinear contact laws generate complex dynamics and persistent difficulties in modelling dissipation. Historically, damping models across these disciplines have been largely empirical \cite{brilliantov1996, antypov2011}, failing to provide consistent restitution control across arbitrary geometries and varying impact conditions.

A rigorously consistent approach to contact mechanics requires that interaction forces derive from scalar potentials defined on the system's configuration space—a core tenet of geometric and variational contact formulations \cite{feng_prep}. This formulation guarantees path-independent work and exact energy conservation for elastic collisions. Recently, an exact phase-space analytical solution was established for the specific case of the power-law damped contact oscillator, demonstrating that the nonlinear trajectory could be mapped to a linear spring-dashpot system \cite{feng_sub, schwager2007}.

This raises a fundamental question: \textit{Is nonlinear contact dynamics intrinsically nonlinear, or is this nonlinearity merely a consequence of the chosen spatial representation?}

We answer this by establishing a hidden harmonic structure that applies not just to power laws, but to all monotone one-dimensional contact systems within this variational framework. While the existence of action–angle variables for one-dimensional conservative systems is a classical result \cite{arnold1989}, the exact harmonic regularisation, universal damping law, and resulting stability bounds derived here are fundamentally new. The key implication is that nonlinear contact dynamics are not fundamentally nonlinear objects; in their natural energy coordinates, they possess an exact harmonic structure. By exploiting this topology, we derive a unique, universal damping law and a rigorous computational stability bound that apply to any arbitrary particle geometry. These results establish a unified structure connecting geometry, dynamics, dissipation, and numerical stability within a single variational framework.

\section{2. Framework Summary: The Main Theorem}

To provide a roadmap for the subsequent derivations, the core results of this paper can be unified into a single structural theorem. The individual results are proved as Theorems 1–3 in Sections 4, 5, and 8.

Let 
\[
H(q,p) = \frac{p^{2}}{2m} + U(q)
\]
be a conservative contact Hamiltonian with $U(0)=0$ and $U'(q)>0$ for $q>0$. 

\begin{enumerate}
    \item The system admits a canonical action–angle representation $(\theta,J)$ in physical time with $\dot{J}=0$ and $\dot{\theta}=dH/dJ$ (Theorem 1).
    \item For any arbitrary positive reference constants $K$ and $M$, it simultaneously admits an exact energy‑based harmonic regularisation $x=\sqrt{2U/K}$ with reparametrised time 
    \[
    \frac{d\tau}{dt}=\sqrt{\frac{M}{m}}\frac{dx}{dq},
    \]
    under which the motion is exactly linear: 
    \[
    Mx''+Kx=0
    \]
    (Theorem 2).
    \item For the dissipative extension 
    \[
    m\ddot{q} + C(q)\dot{q} + U'(q)=0,
    \]
    the transformed dynamics become a linear spring–dashpot 
    \[
    Mx''+C_{0}x'+Kx=0
    \]
    \textbf{if and only if} the physical damping follows the unique law:
    \[
    C(q) \propto U'(q)/\sqrt{U(q)}
    \]
    (Theorem 3).
\end{enumerate}

The proportional form is stated here for clarity; the exact expression with scaling constants is derived in Section 8.

\section{3. Conservative Contact System}

Consider the one-degree-of-freedom conservative contact oscillator  
\[
m\ddot{q} + U'(q)=0,\qquad q\ge 0,
\]
with Hamiltonian 
\begin{equation}
H(q,p) = \frac{p^{2}}{2m} + U(q),
\end{equation}
where $p=m\dot{q}$ and $U(q)$ is the contact potential.

We assume the natural monotonicity conditions for contact mechanics:
(1) $U\in C^{2}((0,\infty))\cap C^{1}([0,\infty))$, 
(2) $U(0)=0$, 
(3) $U'(q)>0$ for all $q>0$, and 
(4) For each energy $E>0$, there exists a unique turning point $q_{\max}(E)>0$ such that $U(q_{\max}(E))=E$.

\section{4. Canonical Action–Angle Structure}

\textbf{Theorem 1 (Action variable and period).} Define the action variable  
\begin{equation}
J(E)=\frac{1}{2\pi}\oint p\,dq = \frac{\sqrt{2m}}{\pi}\int_{0}^{q_{\max}(E)}\sqrt{E-U(q)}\,dq.
\end{equation}
Then $J(E)$ is strictly increasing in $E$, and the period of the conservative contact oscillation is 
\begin{equation}
T(E)=2\pi\,\frac{dJ}{dE}.
\end{equation}

\textbf{Proof.} Differentiating $J(E)$ with respect to $E$ via Leibniz’ rule yields:  
\[
\frac{dJ}{dE} = \frac{\sqrt{2m}}{2\pi}\int_{0}^{q_{\max}(E)}\frac{dq}{\sqrt{E-U(q)}}.
\]
The endpoint term vanishes because $E-U(q_{\max})=0$. Using energy conservation $\dot{q}=\sqrt{2(E-U(q))/m}$, the quarter‑period is:
\[
\frac{T(E)}{4} = \int_{0}^{q_{\max}(E)}\frac{dq}{\dot{q}} = \sqrt{\frac{m}{2}}\int_{0}^{q_{\max}(E)}\frac{dq}{\sqrt{E-U(q)}}.
\]
Comparing the expressions gives 
\[
T(E)=2\pi\,\frac{dJ}{dE}.
\]
\hfill $\blacksquare$

Because $J(E)$ is strictly increasing, it can be inverted to define $E=H(J)$. The system therefore admits canonical action–angle variables $(\theta,J)$ on each periodic annulus such that $H=H(J)$, $\dot{J}=0$, and $\dot{\theta}=\Omega(J)=2\pi/T(E)$, in accordance with the classical Arnold–Liouville theorem \cite{arnold1989}.

\section{5. Exact Energy‑Based Harmonic Regularisation}

\textbf{Theorem 2 (Harmonic Regularisation).} Fix arbitrary constants $K>0$ and $M>0$. Define the energy coordinate  
\begin{equation}
x(q)=\sqrt{\frac{2U(q)}{K}},
\end{equation}
and define the reparametrised time $\tau$ by  
\begin{equation}
\frac{d\tau}{dt} = \sqrt{\frac{M}{m}}\,\frac{dx}{dq}.
\end{equation}
Along every nontrivial conservative trajectory, the transformed motion satisfies  
\begin{equation}
M\frac{d^{2}x}{d\tau^{2}}+Kx=0.
\end{equation}

\textbf{Proof.} By definition, 
\[
U(q)=\frac{1}{2} Kx^{2}.
\]
Kinematically, 
\[
\dot{q} = \frac{dq}{dx}\frac{dx}{d\tau}\frac{d\tau}{dt} = \sqrt{\frac{M}{m}}\,\frac{dx}{d\tau}.
\]
Substituting into the energy equation gives:  
\[
E = \frac{1}{2} M\left(\frac{dx}{d\tau}\right)^{2}+\frac{1}{2} Kx^{2}.
\]
Since $E$ is constant in the conservative system, differentiating with respect to $\tau$ yields 
\[
x'(Mx''+Kx)=0,
\]
confirming the exactly harmonic trajectory. \hfill $\blacksquare$

This demonstrates that the apparent nonlinearity of contact dynamics arises entirely from the choice of physical coordinates; in the energy coordinate with reparametrised time, the system is exactly harmonic.

\textbf{Remark.} The constants $K$ and $M$ are arbitrary; they cancel in all physically meaningful expressions derived from the regularisation (e.g., the damping law and the stability bound). This reflects the fact that the harmonic structure is intrinsic and independent of the choice of reference scales.

\section{6. Rigorous Lower Bound for the Critical Timestep}

The harmonic regularisation translates directly into a computational advantage, providing a closed‑form admissible timestep bound for explicit numerical integration without requiring empirical tuning.

In the regularised virtual space $(x, \tau)$, the system is a linear harmonic oscillator with constant natural frequency 
\[
\Omega_0 = \sqrt{K/M}.
\]
The absolute stability limit for standard explicit central‑difference integrators (such as Velocity Verlet) in this linear space is a universal constant:
\[
\Delta\tau_{\text{crit}} = \frac{2}{\Omega_0}.
\]

Mapping this virtual stability limit back to the physical time domain requires the gradient $\frac{d\tau}{dt}$. A physical timestep $\Delta t$ corresponds to a virtual timestep $\Delta\tau = \frac{d\tau}{dt}\,\Delta t$ at each point along the trajectory. To guarantee stability for the entire nonlinear collision, we require that the virtual step never exceeds $\Delta\tau_{\text{crit}}$. The most restrictive condition occurs where $\frac{d\tau}{dt}$ is largest, because that gives the largest virtual step for a given $\Delta t$. Therefore a sufficient condition for stability is
\[
\Delta t \le \frac{\Delta\tau_{\text{crit}}}{\max_{q \in [0, q_{\max}]} (d\tau/dt)}.
\]

Let $\Delta t_{\text{safe}}$ denote the right‑hand side of this inequality. Since any smaller $\Delta t$ will also be safe, $\Delta t_{\text{safe}}$ is a \textbf{lower bound} on the true maximum stable timestep $\Delta t_{\text{crit}}$ of the nonlinear system, i.e., a conservative admissible timestep guaranteeing stability. In other words,
\[
\Delta t_{\text{crit}} \ge \Delta t_{\text{safe}}.
\]

Substituting 
\[
\frac{d\tau}{dt} = \sqrt{\frac{M}{mK}} \frac{U'(q)}{\sqrt{2U(q)}},
\]
the arbitrary scaling constants $M$ and $K$ cancel with $\Omega_0$, yielding the general expression
\begin{equation}
\Delta t_{\text{safe}} = 2\sqrt{m} \left( \max_{q \in [0, q_{\max}]} \frac{U'(q)}{\sqrt{2U(q)}} \right)^{-1}.
\end{equation}

The location of this maximum is dictated entirely by the geometric shape of the contact:

\begin{enumerate}
    \item \textbf{Monotone Stiffening Contacts:} If the geometry satisfies $2U(q)U''(q) \ge [U'(q)]^2$ (which includes all power‑law forces with exponent $p \ge 1$), the maximum strictly occurs at maximum compression $q_{\max}$. Let the initial impact speed be $v_0$, so that the initial kinetic energy is $E_k = \frac12 m v_0^2$. Energy conservation gives $U(q_{\max}) = E_k$. The bound simplifies to a purely physical expression:
    \begin{equation}
    \Delta t_{\text{safe}} = 2m \frac{v_0}{U'(q_{\max})}.
    \end{equation}
    This provides a closed‑form, sufficient lower bound for the critical timestep.
    \item \textbf{Softening or Purely Linear Contacts:} If the geometry is strictly softening or yields a constant effective stiffness (e.g., a pure volumetric penalty), the maximum gradient occurs at initial contact, requiring the bound to be evaluated at $q \to 0$.
    \item \textbf{Complex Varying Geometries:} For arbitrary shapes that both stiffen and soften, the maximum occurs at an interior point and the supremum must be evaluated explicitly over the collision path.
\end{enumerate}

The derived bound is rigorous: any physical timestep $\Delta t \le \Delta t_{\text{safe}}$ guarantees stability for the nonlinear contact simulation. In practice, $\Delta t_{\text{safe}}$ is a safe upper bound on the timestep (i.e., a lower bound on the true critical timestep) that can be computed directly from the contact geometry and impact conditions.

\section{7. Relation to Canonical Transformations and Action–Angle Variables}

While the regularisation produces a perfectly harmonic oscillator, it is distinct from an ordinary canonical action-angle transformation in standard phase space. 

Applying a point-canonical lift to the spatial transformation 
\[
x = f(q) = \sqrt{2U(q)/K}
\]
yields the conjugate momentum 
\[
P = p/f'(q).
\]
In these coordinates, the Hamiltonian becomes 
\begin{equation}
H(x,P) = \frac{P^{2}}{2M_{\mathrm{eff}}(x)} + \frac{1}{2}Kx^{2},
\end{equation}
where the effective mass is strictly position-dependent: 
\begin{equation}
M_{\mathrm{eff}}(x) = \frac{2mK\,U(q(x))}{[U'(q(x))]^{2}}.
\end{equation}

This reduces to a constant-mass harmonic oscillator in physical time if and only if 
\[
U(q) \propto q^{2}.
\]
For all other nonlinear potentials, the exact harmonic structure appears only after the time reparametrisation $d\tau/dt$. 

\section{8. Dissipative Extension and the Unique Universal Damping Law}

Now consider the dissipative system 
\[
m\ddot{q} + C(q)\dot{q} + U'(q)=0.
\]
Under the transformation $(q,t)\mapsto(x,\tau)$, the system becomes  
\[
Mx''+C_{*}(x)x'+Kx=0,
\]
where 
\[
C_{*}(x)=C(q)\sqrt{\frac{M}{m}}\frac{dq}{dx}.
\]

\textbf{Theorem 3 (Unique universal damping law, necessary and sufficient).} The transformed system is an exact linear spring–dashpot oscillator ($C_{*}(x) \equiv C_{0} > 0$) \textbf{if and only if} the physical damping coefficient satisfies:  
\begin{equation}
C(q)=C_{0}\sqrt{\frac{m}{M}}\frac{U'(q)}{\sqrt{2K\,U(q)}}.
\end{equation}

\textbf{Proof.} By definition, $C_{*}(x)=C(q)\sqrt{M/m}\,dq/dx$. Setting $C_{*}(x)=C_{0}$ gives $C(q)=C_{0}\sqrt{m/M}\,dx/dq$. Substituting $dx/dq = U'(q)/\sqrt{2K U(q)}$ yields the stated condition. \hfill $\blacksquare$ 

This law is the unique dissipation pathway that preserves linear, exponential decay in the regularised harmonic coordinates. It elevates contact damping from an empirical formulation to a rigorous structural requirement for velocity-independent restitution.

\section{9. Recovery of the Power‑Law Family}

For the standard power-law elastic family 
\[
U(q)=\frac{k}{p+1}q^{p+1},
\]
the force scales as $U'(q) = kq^p$. Evaluating the required geometric ratio yields:
\[
\frac{U'(q)}{\sqrt{U(q)}} \propto q^{(p-1)/2}.
\]

The universal damping law therefore demands 
\begin{equation}
C(q)\propto q^{(p-1)/2}.
\end{equation}
This recovers the Tsuji‑type damping exponent \cite{tsuji1992}, confirming that the recently reported exact phase-space solution for power-law damped contacts \cite{feng_sub} is the exact monomial realisation of this universal theorem.

\section{10. Exact Implementation for Arbitrary Geometries}

To demonstrate that the theory seamlessly handles arbitrary, non-power-law geometries at finite penetrations, we first establish the mechanics for a general volumetric shape. Under an $\alpha$-regularised volumetric contact model, the potential energy $U(\delta)$ at penetration depth $\delta$ is defined as:
\begin{equation}
U(\delta) = \frac{K_n}{\alpha+1} V_c(\delta)^{\alpha+1},
\end{equation}
where $V_c(\delta)$ is the exact overlap volume, $K_n$ is the stiffness parameter, and $\alpha$ is a regularisation exponent (e.g., $\alpha = 0.5$ yields quadratic scaling). The elastic restoring force is strictly derived from this potential as $U'(\delta) = K_n V_c(\delta)^\alpha S_n(\delta)$, where $S_n(\delta) = dV_c/d\delta$ is the exact geometric cross-sectional contact area.

To evaluate this without small-overlap asymptotic approximations, consider an ellipsoid with semi-axes \textbf{$a = 0.015$~m} and \textbf{$b = c = 0.008$~m}, and mass \textbf{$m = 0.05$~kg}, impacting a flat wall along the $a$-axis. The exact closed-form expressions for the finite geometric overlap are:
\begin{equation}
S_n(\delta) = \frac{\pi bc}{a^2}(2a\delta - \delta^2), \quad V_c(\delta) = \frac{\pi bc}{a^2}\left(a\delta^2 - \frac{\delta^3}{3}\right).
\end{equation}

Substituting this exact geometry into the universal damping law yields the explicit damping requirement for the ellipsoid:
\begin{equation}
C(\delta) \propto \frac{U'(\delta)}{\sqrt{U(\delta)}} \propto \frac{V_c(\delta)^\alpha S_n(\delta)}{V_c(\delta)^{(\alpha+1)/2}} = S_n(\delta) V_c(\delta)^{(\alpha-1)/2}.
\end{equation}

For the recommended quadratic scaling ($\alpha = 0.5$), this yields the exact, closed-form damping coefficient $C(\delta) \propto S_n(\delta) V_c(\delta)^{-1/4}$. This is completely regular, requires no small-overlap approximations, and maps the complex polynomial geometry of the ellipsoid impact exactly onto a linear harmonic oscillator under the energy-based transformation.

\begin{figure}[t]
\includegraphics[width=\columnwidth]{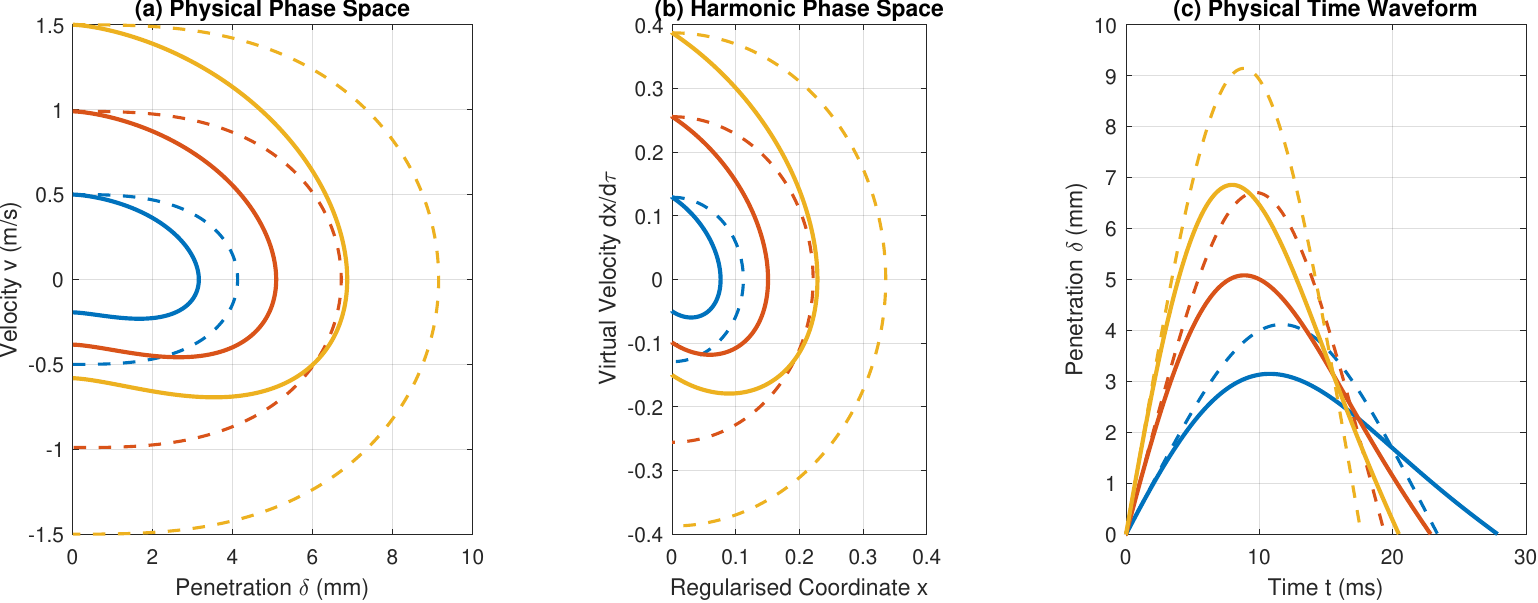}
\caption{\label{fig:fig1}(Color online) Exact harmonic regularisation of a highly nonlinear contact oscillator (an ellipsoid impacting a flat wall with $\alpha = 0.5$ volumetric regularisation and stiffness \textbf{$K_n = 10^8$~Pa}). Dashed lines denote the strictly conservative interaction ($C_0 = 0$), while solid lines denote the dissipative extension governed by the universal damping law ($C_0 = 0.5$). The colours correspond to the initial impact velocities: \textbf{$v_0 = 0.50$~m/s} (blue), \textbf{$0.99$~m/s} (red), and \textbf{$1.50$~m/s} (yellow). \textbf{(a)} In the physical phase space, the conservative trajectories are symmetric but non-elliptical, whereas the dissipative trajectories are heavily skewed and asymmetric. \textbf{(b)} Under the exact energy-based coordinate transformation $x(\delta)$ and time reparametrisation $\tau$ (using virtual scaling constants $K=1$ and $M=0.75$), the identical physical trajectories map perfectly onto the symmetric semi-ellipses (conservative) and exact logarithmic spirals (dissipative) of a linear harmonic oscillator. \textbf{(c)} The physical time-domain waveforms illustrate the corresponding nonlinear penetration pulses and energy loss.}
\end{figure}

The numerical validation of this exact mapping is presented in Fig.~\ref{fig:fig1} for three distinct impact velocities (\textbf{0.50~m/s}, \textbf{0.99~m/s}, and \textbf{1.50~m/s}) with \textbf{$K_n = 10^8$~Pa}. In the physical phase space (Fig.~\ref{fig:fig1}a), the conservative trajectories are symmetric but distinctly non-elliptical, while the dissipative trajectories exhibit heavily skewed, asymmetric loops. However, when mapped into the regularised harmonic phase space (Fig.~\ref{fig:fig1}b), the conservative trajectories (dashed lines) form perfect, symmetric semi-ellipses, validating Theorem 2. Simultaneously, the dissipative trajectories governed by the universal damping law (solid lines) transform into mathematically exact logarithmic spirals, visually confirming that the framework perfectly linearises the dissipation (Theorem 3). Numerical verification, performed using MATLAB's \textsf{ode45} integrator with stringent tolerances (\textsf{RelTol} $= 10^{-13}$, \textsf{AbsTol} $= 10^{-13}$), confirms that the Hamiltonian energy in this transformed space is conserved to machine precision (relative error $< 3 \times 10^{-13}$), guaranteeing the exactness of the harmonic structure independent of graphical representation.

Furthermore, this $\alpha$-regularised ellipsoid model strictly satisfies the stiffening condition $2U(\delta)U''(\delta) \ge [U'(\delta)]^2$ across the entire physically valid elastic collision domain (up to a penetration of $\delta \le 0.735a$). Consequently, the maximum transformation gradient is guaranteed to occur at maximum compression, allowing the rigorous critical timestep to be computed directly using the exact kinetic energy shortcut derived in Section 6. The resulting absolute lower bounds for numerical stability are summarised in Table~\ref{tab:table1}.

\begin{table}[b]
\caption{\label{tab:table1} Theoretical lower bounds for the critical timestep ($\Delta t_{\text{safe}}$) evaluated at maximum compression for the ellipsoid impact scenarios.}
\begin{ruledtabular}
\begin{tabular}{lccc}
Impact Velocity & Max Penetration & Max Force & Stability Bound \\
$v_0$ (m/s) & $\delta_{\max}$ (mm) & $U'(\delta_{\max})$ (N) & $\Delta t_{\text{safe}}$ (ms) \\
\colrule
0.50 & 4.116586 & 4.325435 & 11.559531 \\
0.99 & 6.706989 & 10.000023 & 9.899977 \\
1.50 & 9.143383 & 16.102848 & 9.315123 \\
\end{tabular}
\end{ruledtabular}
\end{table}

\section{11. Conclusion}

We have demonstrated that the widely accepted nonlinearity of conservative and dissipative contact dynamics is not an intrinsic physical property, but a consequence of the spatial coordinate representation. By elevating the problem to an energy-based phase space with a reparametrised time, any monotone energy-consistent contact model resolves into an exact linear harmonic oscillator.

This regularisation yields two immediate and powerful consequences for granular mechanics and DEM simulations. First, it dictates a unique universal damping law that preserves this linear topology, providing exact velocity-independent restitution across arbitrary particle geometries while recovering classical empirical models as monomial special cases. Second, it yields a rigorous closed-form sufficient bound for numerical stability, eliminating the need for empirical timestep tuning.

Beyond granular simulations, this exact regularisation provides a unified analytical foundation for dissipative acoustic metamaterials, viscoelastic AFM indentation, and nanoscale friction modelling. Ultimately, this framework establishes nonlinear contact dynamics as a class of systems admitting an exact linear representation in a transformed phase space, thereby embedding dissipative contact mechanics within the structure-preserving framework of analytical mechanics.


\begin{thebibliography}{99}

\bibitem{hertz1882}
H. Hertz, J. Reine Angew. Math. \textbf{92}, 156 (1882).

\bibitem{tsuji1992}
Y. Tsuji, T. Tanaka, and T. Ishida, Powder Technol. \textbf{71}, 239 (1992).

\bibitem{daraio2006}
C. Daraio, V. F. Nesterenko, E. B. Herbold, and S. Jin, Phys. Rev. Lett. \textbf{96}, 058002 (2006).

\bibitem{giessibl2003}
F. J. Giessibl, Rev. Mod. Phys. \textbf{75}, 949 (2003).

\bibitem{garcia2002}
R. Garcia and R. Perez, Surf. Sci. Rep. \textbf{47}, 197 (2002).

\bibitem{brilliantov1996}
N. V. Brilliantov, F. Spahn, J.-M. Hertzsch, and T. P\"{o}schel, Phys. Rev. E \textbf{53}, 5382 (1996).

\bibitem{antypov2011}
D. Antypov and J. A. Elliott, EPL \textbf{94}, 50004 (2011).

\bibitem{feng_prep}
Y. T. Feng, \textit{Geometric Contact Mechanics for Particle Systems} (in preparation).

\bibitem{feng_sub}
Y. T. Feng, https://arxiv.org/abs/2603.27764 (2026).

\bibitem{schwager2007}
T. Schwager and T. P\"{o}schel, Granular Matter \textbf{9}, 465 (2007).

\bibitem{arnold1989}
V. I. Arnold, \textit{Mathematical Methods of Classical Mechanics} (Springer, New York, 1989).

\end{thebibliography}
\end{document}